\documentclass{mfatshort}
\usepackage{cite}

\newtheorem{theorem}{Theorem}[section]

\newtheorem{proposition}{Proposition}[section]

\theoremstyle{remark}

\theoremstyle{definition}

\newtheorem{definition}{Definition}[section]

\numberwithin{equation}{section}

\begin{document}
\title[A condition for solutions of a parabolic
problem to be classical]{A condition for generalized solutions\\ of a parabolic problem for a Petrovskii system\\ to be classical}

\author[Valerii Los]{Valerii Los}

\address{National Technical University of Ukraine Igor Sikorsky Kyiv Polytechnic Institute,
Prospect Peremohy 37, 03056, Kyiv-56, Ukraine}

\email{v$\underline{\phantom{k}}$\,los@yahoo.com}

\subjclass[2010]{35K35, 46E35}


\dedicatory{Dedicated to the 70th anniversary of Professor V. A. Mikhailets}

\keywords{Parabolic problem, H\"ormander space, slowly varying function, generalized solution, classical solution.}

\begin{abstract}
We obtain a new sufficient condition under which generalized solutions to a parabolic initial-boundary-value problem for a Petrovskii system and the homogeneous Cauchy data are classical. The condition is formulated in terms of the belonging of the right-hand sides of the problem to some anisotropic H\"ormander spaces.
\end{abstract}

\maketitle

\section{Introduction}\label{sec1}
In the theory of partial differential equations, of great importance are explicit conditions that guarantee a required regularity of solutions to the equations under study. As a rule, initial-boundary-value problems are investigated on appropriate pairs of normed distribution spaces. Once a solvability result for a problem is obtained, it is natural to ask when a distributional (i.e. generalized) solution to this problem is classical.  In other words, when one may calculate the left-hand sides of the problem via classical partial derivatives and via traces of continuous functions. An answer to this question is usually given in terms of the belonging of the right-hand sides of the problem to relevant distribution spaces. The more finely a scale of these spaces is calibrated, the more exact conditions will be got.

Investigating parabolic problems, one usually uses anisotropic Sobolev spaces or H\"older spaces parametrized with a pair of real numbers \cite{LionsMagenes72ii, Ivasyshen90, Lunardi1995, ZhitarashuEidelman98,Solonnikov65}. To obtain more precise results about regularity of solutions of these problems, it is natural to resort to  distribution spaces characterized with the help of function parameters. Broad classes of such spaces were introduced and investigated by H\"ormander \cite[Section~2.2]{Hermander63} and Volevich and Paneah \cite{VolevichPaneah65}. Of late years a theory of solvability of general parabolic problems is developed for some Hilbert anisotropic H\"ormander spaces \cite{Los16JMathSci, Los16UMJ6, Los17UMJ3, LosMikhailetsMurach17CPAA, LosMurach13MFAT2, LosMurach17OpenMath, Los17MFAT2, LosMikhailetsMurach19arXiv}. Their order of regularity is given by a pair of real numbers $s$ and $s/(2b)$, where $2b$ is the parabolic weight of the problem, and by a function $\varphi:[1,\infty)\rightarrow(0,\infty)$ that varies regularly at infinity. The function parameter $\varphi$ defines additional (positive or negative) regularity of distributions forming these spaces. The $\varphi(\cdot)\equiv1$ case provides the anisotropic Sobolev space of order $s$ with respect to the spatial variables and order $s/(2b)$ with respect to the time variable. A core of this theory consists of isomorphism theorems for operators induced by parabolic problems and acting between appropriate H\"ormander spaces. (Somewhat earlier a theory of elliptic boundary-value problems was built for isotropic versions of these spaces \cite{MikhailetsMurach12BJMA2, MikhailetsMurach14}.)
\\

The present paper investigates a parabolic initial-boundary-value problem for a Petrovskii parabolic system and the homogeneous Cauchy data. An isomorphism theorem for operators generated by this problem on pairs of appropriate inner product H\"ormander spaces is proved in \cite{Los17UMJ3}. The purpose of the present paper is to supplement this theorem with a corresponding sufficient condition for generalized solutions of the parabolic problem to be classical in the sense mentioned above. The condition is formulated in terms of the belonging of the right-hand sides of the problem to suitable anisotropic H\"ormander spaces. The use of the function parameter $\varphi$ allows us to achieve the minimal admissible value of the number parameter $s$, which is not possible in the framework of Sobolev spaces or H\"older spaces \cite{Ilin60, IlinKalashOleinik62, Solonnikov65}. As to scalar  parabolic problems, conditions of this type are obtained in \cite{Los16UMJ9, Los16UMJ11}.

\section{Statement of the problem}\label{sec2}
We arbitrarily choose an integer $n\geq2$ and a real number $\tau>0$. Let $G$ be a bounded domain in $\mathbb{R}^{n}$ with an infinitely smooth boundary $\Gamma:=\partial G$. We put
$\Omega:=G\times(0,\tau)$ and $S:=\Gamma\times(0,\tau)$; thus, $\Omega$ is an open cylinder in $\mathbb{R}^{n+1}$, and $S$ is its lateral boundary. Then $\overline{\Omega}:=\overline{G}\times[0,\tau]$ and
$\overline{S}:=\Gamma\times[0,\tau]$ are the closures of $\Omega$ and $S$ respectively.

We consider the following parabolic initial-boundary-value problem in $\Omega$:
\begin{equation}\label{24f1}
\begin{split}
\sum_{k=1}^{N}A_{j,k}(x,t,D_x,&\partial_t)u_{k}(x,t)=f_{j}(x,t)\\
&\mbox{for all}\quad x\in G,\quad t\in(0,\tau)\quad\mbox{and}\quad
j\in\{1,\dots,N\};
\end{split}
\end{equation}
\begin{equation}\label{24f2}
\begin{split}
\sum_{k=1}^{N}B_{j,k}(x,t,D_x,&\partial_t)u_{k}(x,t)\big|_{S}=
g_{j}(x,t)\\
&\mbox{for all}\quad x\in\Gamma,\quad 0<t<\tau \quad\mbox{and}\quad
j\in\{1,\dots,m\};
\end{split}
\end{equation}
\begin{equation}\label{24f3}
\begin{split}
\partial^{\,r}_t u_{k}(&x,t)\big|_{t=0}=0\\
&\mbox{for all}\quad x\in G,\quad
k\in\{1,\ldots,N\}\quad\mbox{and}\quad r\in\{0,\ldots,\varkappa_k-1\}.
\end{split}
\end{equation}
Note that the initial data \eqref{24f3} are assumed to be zero.
The linear partial differential operators (PDOs) used in the problem are of the form
\begin{gather}\label{24f4}
A_{j,k}(x,t,D_x,\partial_t):=
\sum_{|\alpha|+2b\beta\leq 2b\varkappa_{k}}
a^{\alpha,\beta}_{j,k}(x,t)\,D^\alpha_x\partial^\beta_t,\\
B_{j,k}(x,t,D_x,\partial_t):=
\left\{
\begin{array}{cl}
\sum\limits_{|\alpha|+2b\beta\leq l_j+2b\varkappa_{k}}
b^{\alpha,\beta}_{j,k}(x,t)\,D^\alpha_x\partial^\beta_t& \mbox{if}\quad l_j+2b\varkappa_{k}\geq0,\\
0 & \mbox{if}\quad l_j+2b\varkappa_{k}<0
\end{array}
\right. \label{24f5}
\end{gather}
for all admissible values of the indexes $j$ and $k$. Here, the positive integers $N\geq2$, $b$, and $\varkappa_{1},\ldots,\varkappa_{N}$ are arbitrarily chosen; $m:=b(\varkappa_{1}+\dots+\varkappa_N)$, and  $l_1,\ldots,l_m\in\mathbb{Z}$. The even number $2b$ is called the parabolic weight of this problem. All coefficients of the PDOs $A_{j,k}$ and $B_{j,k}$ are supposed to be infinitely smooth complex-valued functions given on $\overline{\Omega}$ and $\overline{S}$ respectively; i.e., each
\begin{equation*}
a_{j,k}^{\alpha,\beta}\in C^{\infty}(\overline{\Omega}):=
\bigl\{w\!\upharpoonright\overline{\Omega}\!:\,w\in C^{\infty}(\mathbb{R}^{n+1})\bigr\}
\end{equation*}
and each
\begin{equation*}
b_{j,k}^{\alpha,\beta}\in C^{\infty}(\overline{S}):=
\bigl\{v\!\upharpoonright\overline{S}\!:\,v\in C^{\infty}(\Gamma\times\mathbb{R})\bigr\}.
\end{equation*}

We use the notation
$D^\alpha_x:=D^{\alpha_1}_{1}\dots D^{\alpha_n}_{n}$, with $D_{k}:=i\,\partial/\partial{x_k}$, and $\partial_t:=\partial/\partial t$
for the partial derivatives of functions depending on $x=(x_1,\ldots,x_n)\in\mathbb{R}^{n}$ and $t\in\mathbb{R}$. Here, $i$ is imaginary unit, and $\alpha=(\alpha_1,...,\alpha_n)$ is a multi-index, with $|\alpha|:=\alpha_1+\cdots+\alpha_n$. In formulas \eqref{24f4} and \eqref{24f5} and their analogs, we take summation over the integer-valued nonnegative indices $\alpha_1,...,\alpha_n$ and $\beta$ that satisfy the condition written under the integral sign.

We assume that the initial-boundary value problem \eqref{24f1}--\eqref{24f3} is Petrovskii parabolic in the cylinder $\Omega$. Let us recall the corresponding definition \cite[Section~1, \S~1]{Solonnikov65}. Define the principal symbols of the PDOs \eqref{24f4} and \eqref{24f5} as follows:
\begin{equation*}
A^{(0)}_{j,k}(x,t,\xi,p):=\sum_{|\alpha|+2b\beta=2b\varkappa_{k}}
a^{\alpha,\beta}_{j,k}(x,t)\,\xi^\alpha p^\beta,
\end{equation*}
\begin{equation*}
B^{(0)}_{j,k}(x,t,\xi,p)=:
\left\{
\begin{array}{cl}
\sum\limits_{|\alpha|+2b\beta=l_j+2b\varkappa_{k}}
b^{\alpha,\beta}_{j,k}(x,t)\,\xi^\alpha p^\beta,& \mbox{if}\quad l_j+2b\varkappa_{k}\geq 0,
\\ 0, & \mbox{if}\quad l_j+2b\varkappa_{k}<0.
\end{array}
\right.
\end{equation*}
These symbols are homogeneous polynomials in
$\xi:=(\xi_{1},\ldots,\xi_{n})\in\mathbb{C}^{n}$ and $p\in\mathbb{C}$ jointly (as usual, $\xi^\alpha:=\xi_{1}^{\alpha_{1}}\ldots\xi_{n}^{\alpha_{n}}$).
Consider the matrices
\begin{equation*}
A^{(0)}(x,t,\xi,p):=\bigl(A^{(0)}_{j,k}(x,t,\xi,p)\bigr)_{j,k=1}^{N},
\end{equation*}
\begin{equation*}
B^{(0)}(x,t,\xi,p):=\bigl(B^{(0)}_{j,k}(x,t,\xi,p)\bigr)
_{\substack{j=1,\ldots,m\\k=1,\ldots,N}}.
\end{equation*}

The problem \eqref{24f1}--\eqref{24f3} is said to be Petrovskii parabolic in $\Omega$ if it satisfies the following three conditions:
\begin{itemize}
\item [(i)] For arbitrary points $x\in\overline{G}$ and $t\in[0,\tau]$ and every vector $\xi\in\mathbb{R}^{n}$, all the roots $p(x,t,\xi)$ of the polynomial $\det A^{(0)}(x,t,\xi,p)$ in $p\in\mathbb{C}$ satisfy the inequality $\mathrm{Re}\,p(x,t,\xi)\leq -\delta\,|\xi|^{2b}$ for some number $\delta>0$ that does not depend on $x$, $t$, and~$\xi$.
\item [(ii)] Each equation in the system \eqref{24f1} is solvable with respect to the derivative $\partial^{\varkappa_j}_{t}u_{j}$, where $j$ is the number of this equation, and does not contain any derivative of the form $\partial^{\varkappa_k}_{t}u_{k}$ where $k\neq j$. Thus, we may and do assume that
$a_{j,k}^{(0,0,\dots,0),\varkappa_{k}}(x,t)\equiv\delta_{j,k}$ whenever $j,k\in\{1,\dots,N\}$ (as usual, $\delta_{j,k}$ is the Kronecker delta).
\end{itemize}

To formulate the third condition, we fix a number $\delta_1\in(0,\delta)$, where $\delta$ has appeared in Condition~(i), and then arbitrarily choose a point $x\in\Gamma$, real number $t\in[0,\tau]$, vector $\xi\in\mathbb{R}^{n}$ tangent to the boundary $\Gamma$ at $x$, and number $p\in\mathbb{C}$ such that $\mathrm{Re}\,p\geq -\delta_1|\xi|^{2b}$ and $|\xi|+|p|\neq0$. Let $\nu(x)$ denote the unit vector of the inward normal to $\Gamma$ at $x$. It follows from Condition~(i) and the inequality $n\geq2$ that the polynomial
$\det A^{(0)}(x,t,\xi+\zeta\nu(x),p)$ in $\zeta\in\mathbb{C}$ has $m$ roots $\zeta^{+}_{j}(x,t,\xi,p)$, $j=\nobreak1,\ldots,m$, with positive imaginary part and $m$ roots with negative imaginary part provided that each root is taken the number of times equal to its multiplicity.

The third condition is formulated as follows:
\begin{itemize}
\item [(iii)] For each positive number $\delta_1<\delta$ and for every choice of the parameters $x$, $t$, $\xi$ ³ $p$ indicated above, the rows of the matrix
\begin{equation*}
B^{(0)}(x,t,\xi+\zeta\nu(x),p)\cdot\widetilde{A}^{(0)}
(x,t,\xi+\zeta\nu(x),p)
\end{equation*}
are linearly independent modulo the polynomial
$\prod_{j=1}^{m}(\zeta-\zeta^{+}_{j}(x,t,\xi,p))$. Here,  $\widetilde{A}^{(0)}$ is the transposed matrix of the cofactors of entries of $A^{(0)}$.
\end{itemize}

Note that Conditions~(i) and (ii) means that the system \eqref{24f1} is uniformly $2b$-parabolic in the sense of Petrovskii in $\overline{\Omega}$  \cite{Petrovskii38}, whereas Condition~(iii) claims that the collection of boundary conditions \eqref{24f2} covers the parabolic system \eqref{24f1} on $\overline{S}$.

\section{H\"ormander spaces related to the parabolic problem}\label{sec3}

Following \cite[Section~3]{LosMikhailetsMurach17CPAA}, we consider H\"ormander spaces used for the investigation of the parabolic problem \eqref{24f1}--\eqref{24f3}. They are parametrized with two numbers $s$ and $s\gamma$, where $s\in\mathbb{R}$ and $\gamma:=1/(2b)$, and with a function $\varphi\in\mathcal{M}$. The class $\mathcal{M}$ is defined to consist of all Borel measurable functions
$\varphi:[1,\infty)\rightarrow(0,\infty)$ such that
\begin{itemize}
\item [($\ast$)] both the functions $\varphi$ and $1/\varphi$ are bounded on each
compact interval $[1,d]$, with $1<d<\infty$;
\item [($\ast\ast$)] $\varphi$ is a slowly varying function at infinity in the sense of J.~Karamata \cite{Karamata30a}; i.e.,
    $\varphi(\lambda r)/\varphi(r)\to1$, as $r\to\infty$, whenever $\lambda>0$.
\end{itemize}

By definition, the complex linear space $H^{s,s\gamma;\varphi}(\mathbb{R}^{k+1})$, where $1\leq k\in\mathbb{Z}$, consists of all tempered distributions $w$ on $\mathbb{R}^{k+1}$ whose (complete) Fourier transform $\widetilde{w}$ is locally Lebesgue integrable over $\mathbb{R}^{k+1}$ and satisfies the condition
\begin{equation}\label{norm}
\|w\|_{H^{s,s\gamma;\varphi}(\mathbb{R}^{k+1})}:=
\biggl(\;\int\limits_{\mathbb{R}^{k}}\int\limits_{\mathbb{R}}
r_{\gamma}^{2s}(\xi,\eta)\,\varphi^{2}(r_{\gamma}(\xi,\eta))\,
|\widetilde{w}(\xi,\eta)|^{2}\,d\xi\,d\eta\biggr)^{1/2}<\infty,
\end{equation}
where
$$
r_{\gamma}(\xi,\eta):=\bigl(1+|\xi|^2+|\eta|^{2\gamma}\bigr)^{1/2}
\quad\mbox{for all}\;\;\xi\in\mathbb{R}^{k}\;\;\mbox{and}
\;\;\eta\in\mathbb{R}.
$$
This space is Hilbert and separable with respect to the norm \eqref{norm}.

It is a special case of the spaces $\mathcal{B}_{p,\mu}$ introduced by H\"ormander \cite[Section~2.2]{Hermander63}; namely, $H^{s,s\gamma;\varphi}(\mathbb{R}^{k+1})=
\mathcal{B}_{p,\mu}$ provided that $p=2$ and $\mu(\xi,\eta)\equiv r_{\gamma}^{s}(\xi,\eta)\varphi(r_{\gamma}(\xi,\eta))$. If $\varphi(\cdot)\equiv1$, the space $H^{s,s\gamma;\varphi}(\mathbb{R}^{k+1})$ becomes the anisotropic Sobolev space $H^{s,s\gamma}(\mathbb{R}^{k+1})$. Generally, we have the dense continuous embeddings
\begin{equation}\label{embeddings}
H^{s_{1},s_{1}\gamma}(\mathbb{R}^{k+1})\hookrightarrow
H^{s,s\gamma;\varphi}(\mathbb{R}^{k+1})\hookrightarrow
H^{s_{0},s_{0}\gamma}(\mathbb{R}^{k+1})\quad
\mbox{whenever}\quad s_{0}<s<s_{1}.
\end{equation}

Basing on $H^{s,s\gamma;\varphi}(\mathbb{R}^{k+1})$, consider some Hilbert function spaces  relating to the problem \eqref{24f1}--\eqref{24f3}. Let $V$ be an open nonempty set in $\mathbb{R}^{k+1}$. (Specifically, we need the case where
$V=\Omega$, with $k=n$.) Put
\begin{equation}\label{24f7}
H^{s,s\gamma;\varphi}_{+}(V):=\bigl\{w\!\upharpoonright\!V:\,
w\in H^{s,s\gamma;\varphi}(\mathbb{R}^{k+1}),\;\,\mathrm{supp}\,
w\subseteq\mathbb{R}^{k}\times[0,\infty)\bigl\}.
\end{equation}
The norm in the complex linear space \eqref{24f7} is defined by the formula
\begin{equation}\label{24f8}
\begin{split}
\|u\|_{H^{s,s\gamma;\varphi}_{+}(V)}&:=
\inf\bigl\{\,\|w\|_{H^{s,s\gamma;\varphi}(\mathbb{R}^{k+1})}:\\
&w\in H^{s,s\gamma;\varphi}(\mathbb{R}^{k+1}),\;\,
\mathrm{supp}\,w\subseteq\mathbb{R}^{k}\times[0,\infty),\;\,
u=w\!\upharpoonright\!V\bigl\},
\end{split}
\end{equation}
with $u\in H^{s,s\gamma;\varphi}_{+}(V)$. This space is Hilbert and separable with respect to this norm. Specifically, the set
$$
C^{\infty}_{+}(\overline{\Omega}):=
\bigl\{w\!\upharpoonright\overline{\Omega}:\,
w\in C^{\infty}(\mathbb{R}^{n+1}),\;\,
\mathrm{supp}\,w\subseteq\mathbb{R}^{n}\times[0,\infty)\bigr\}
$$
is dense in $H^{s,s\gamma;\varphi}_{+}(\Omega)$.

We also consider the space $H^{s,s\gamma;\varphi}_{+}(S)$ on the lateral boundary $S$ of the cylinder $\Omega$, we restricting ourselves to the $s>0$ case. Briefly saying, this space consists of all functions $v\in L_2(S)$ that yield functions from the space $H^{s,s\gamma;\varphi}_{+}(\Pi)$ on $\Pi:=\mathbb{R}^{n-1}\times(0,\tau)$ with the help of some local coordinates on $\overline{S}$. Let us turn to a detailed definition.

We arbitrarily choose a finite atlas on $\Gamma$ of class $C^{\infty}$.
Let this atlas be formed by some local charts $\theta_{j}:\mathbb{R}^{n-1}\leftrightarrow\Gamma_{j}$, with
$j=1,\ldots,\lambda$. Here, $\Gamma_{1},\ldots,\Gamma_{\lambda}$ are open nonempty subsets of $\Gamma$ such that $\Gamma:=\Gamma_{1}\cup\cdots\cup\Gamma_{\lambda}$. We also arbitrarily choose functions
$\chi_{j}\in C^{\infty}(\Gamma)$, with $j=1,\ldots,\lambda$, such that
$\mathrm{supp}\,\chi_{j}\subset\Gamma_{j}$ and $\chi_{1}+\cdots+\chi_{\lambda}=1$ on $\Gamma$. Thus, these functions form a partition of unity on $\Gamma$.

By definition, the complex linear space $H^{s,s\gamma;\varphi}_{+}(S)$ consists of all functions $v\in L_2(S)$ such that the function $v_{j}(y,t):=\chi_{j}(\theta_{j}(y))v(\theta_{j}(y),t)$ of $y\in\mathbb{R}^{n-1}$ and $t>0$ belongs to $H^{s,s\gamma;\varphi}_{+}(\Pi)$ for each $j\in\{1,\ldots,\lambda\}$. (As usual, $L_2(S)$ denotes the space of all square integrable functions on the surface $S$.) The space $H^{s,s\gamma;\varphi}_{+}(S)$ is separable Hilbert with respect to the norm
$$
\|v\|_{H^{s,s\gamma;\varphi}_{+}(S)}:=
\bigl(\|v_{1}\|_{H^{s,s\gamma;\varphi}_{+}(\Pi)}^{2}+\cdots+
\|v_{\lambda}\|_{H^{s,s\gamma;\varphi}_{+}(\Pi)}^{2}\bigr)^{1/2}.
$$
This space does not depend up to equivalence of norms on the indicated choice of an atlas and partition of unity on $\Gamma$.

To formulate the main result of the paper, we need local versions of the spaces just considered. Let $U$ be an open subset of $\mathbb{R}^{n+1}$ such that $\Omega_0:=U\cap\Omega\neq\emptyset$, and put $\Omega':=U\cap\partial\Omega$, $S_0:=U\cap S$, and $S':=U\cap \partial S$. We let $H^{s,s\gamma;\varphi}_{+,\mathrm{loc}}(\Omega_0,\Omega')$ denote the linear space of all distributions $u\in\mathcal{D}'(\Omega)$
such that $\chi u\in H^{s,s\gamma;\varphi}_{+}(\Omega)$ for every function $\chi\in C^\infty (\overline\Omega)$ subject to $\mathrm{supp}\,\chi\subset\Omega_0\cup\Omega'$. Analogously, $H^{s,s\gamma;\varphi}_{+,\mathrm{loc}}(S_0,S')$ denotes the linear space of all distributions $v\in\mathcal{D}'(S)$ such that $\chi v\in H^{s,s\gamma;\varphi}_{+}(S)$ for every function $\chi\in C^\infty (\overline S)$ subject to $\mathrm{supp}\,\chi\subset S_0\cup S'$. Here, as usual, $\mathcal{D}'(\Omega)$ and $\mathcal{D}'(S)$ stand for the linear topological spaces of all distributions on $\Omega$ or $S$, resp.

If $\varphi(\cdot)\equiv1$, we will omit the index $\varphi$ in the designations of the spaces considered in this section.

\section{Main result}\label{sec4}
Analyzing the problem \eqref{24f1}--\eqref{24f3}, we put $u:=(u_1,\ldots,u_N)$, $f:=(f_1,\ldots,f_N)$, and $g:=(g_1,\ldots,g_m)$. Let $\sigma_0:=\max\{0,l_1+1,\dots,l_m+1\}$. Consider the mapping $(C^{\infty}_{+}(\overline{\Omega}))^{N}\ni u\mapsto(f,g)$, where $f$ and $g$ are defined by \eqref{24f1} and \eqref{24f2}. It follows from \cite[Theorem~5.7]{ZhitarashuEidelman98}, that this mapping extends uniquely (by continuity) to an isomorphism
\begin{align*}
\Lambda:\:&\mathcal{G}_{+}(\Omega):=
\bigoplus_{k=1}^{N}
H^{\sigma_0+2b\varkappa_k,(\sigma_0+2b\varkappa_k)/(2b)}_{+}(\Omega)\\
&\leftrightarrow\bigl(H^{\sigma_0,\sigma_0/(2b)}_{+}(\Omega)\bigr)^N\oplus
\bigoplus_{j=1}^{m}H^{\sigma_0-l_j-1/2,(\sigma_0-l_j-1/2)/(2b)}_{+}(S).
\end{align*}

Let $f\in(\mathcal{D}'(\Omega))^{N}$ and $g\in(\mathcal{D}'(S))^{m}$; then a vector function $u\in \mathcal{G}_{+}(\Omega)$ is said to be a (strong) generalized solution of the problem \eqref{24f1}--\eqref{24f3} if $\Lambda u=(f,g)$.

To give a notion of a classical solution of this problem, we let $l_0:=\max\{l_1,\dots,l_m\}$ and put $S_{\varepsilon}:=\{x\in\Omega:\mbox{dist}(x,S)<\varepsilon\}$ and
$G_{\varepsilon}:=\{x\in\Omega:\mbox{dist}(x,G)<\varepsilon\}$ for sufficiently small $\varepsilon>0$.

\begin{definition}\label{main-def}
A generalized solution $u\in \mathcal{G}_{+}(\Omega)$ of the problem
\eqref{24f1}--\eqref{24f3} is called \emph{classical} if generalized partial derivatives of each scalar function $u_k=u_k(x,t)$, where $k\in\{1,\dots,N\}$, satisfy the following conditions:
\begin{itemize}
\item [(a)] $D^\alpha_x\partial^\beta_t u_k$ is continuous on $\Omega$ whenever
$0\leq|\alpha|+2b\beta\leq 2b\varkappa_k$;
\item [(b)] $D^\alpha_x\partial^\beta_t u_k$ is continuous on $S_{\varepsilon}\cup S$ for some $\varepsilon>0$ whenever $0\leq|\alpha|+2b\beta\leq l_0+2b\varkappa_k$;
\item [(c)] $\partial^r_t u_k$ is continuous on $G_{\varepsilon}\cup G$ for some $\varepsilon>0$ whenever $0\leq r\leq\varkappa_k-1$.
\end{itemize}
\end{definition}

This definition claims the minimal conditions under which the left-hand sides of the problem are calculated with the help of continuous classical derivatives. Note that conditions (b) and (c) allow us to take the traces on $S$ and $G$ in \eqref{24f2} and \eqref{24f3} in the sense of  restriction of a continuous function.

The main result of the paper reads as follows:

\begin{theorem}\label{24th4.1}
Let $\sigma_1:=b+n/2$, $\sigma_2:=l_0+b+n/2$, and $\sigma_3:=-b+n/2$, and suppose that $\sigma_2>\sigma_0$ and $\sigma_3>\sigma_0$.
Assume that a vector function $u\in\mathcal{G}_{+}(\Omega)$ is a generalized solution of the parabolic problem \eqref{24f1}--\eqref{24f3} whose right-hand sides satisfy the conditions
\begin{equation}\label{24f12}
f\in(H^{\sigma_1,\sigma_1/(2b);\varphi_1}_{+,\mathrm{loc}}
(\Omega,\emptyset))^N\cap
(H^{\sigma_2,\sigma_2/(2b);\varphi_2}_{+,\mathrm{loc}}
(S_{\varepsilon},S))^N\cap (H_{+,\mathrm{loc}}^{\sigma_3,\sigma_3/(2b);\varphi_3}
(G_{\varepsilon},G))^N
\end{equation}
and
\begin{equation}\label{24f12bis}
g\in\bigoplus_{j=1}^{m} H_{+,\mathrm{loc}}^{\sigma_2-l_j-1/2,(\sigma_2-l_j-1/2)/(2b);\varphi_2}
(S,\emptyset)
\end{equation}
for some function parameters
$\varphi_1$, $\varphi_2$, and $\varphi_3$ such that
\begin{equation}\label{24f13}
\int\limits_{1}^{\,\infty}\;\frac{dr}{r\varphi_j^2(r)}<\infty
\quad\mbox{for each}\quad j\in\{1,2,3\}.
\end{equation}
Then $u$ is a classical solution of this problem.
\end{theorem}

The use of H\"ormander spaces allows us to attain the minimal admissible values of the number parameters in conditions \eqref{24f12} and \eqref{24f12bis}. If we formulate an analog of this theorem using anisotropic Sobolev spaces (i.e. restricting ourselves to the case where $\varphi_1=\varphi_2=\varphi_3=1$), we have to claim that the right-hand sides of the problem under investigation satisfy these conditions for certain $\sigma_1>b+n/2$,  $\sigma_2>l_0+b+n/2$, and $\sigma_3>-b+n/2$.

\section{Proof of the main result}\label{sec6}

The proof is based on the following regularity property of the generalized solutions to the considered problem \cite[Theorem~3]{Los17UMJ3}:

\begin{proposition}\label{24prop1}
Let an integer $k\in\{1,\dots,N\}$, and let an integer $p\geq0$ satisfy
\begin{equation}\label{p-condition}
p+b+n/2>\sigma_0+2b\varkappa_k
\end{equation}
Assume that a vector function
$u\in\mathcal{G}_+(\Omega)$ is a generalized solution of the parabolic problem \eqref{24f1}--\eqref{24f3} whose right-hand sides satisfy the conditions
\begin{equation}\label{24f14}
f\in\bigl(H^{\sigma,\sigma/(2b);\varphi}_{+,\mathrm{loc}}
(\Omega_0,\Omega')\bigr)^N
\end{equation}
and
\begin{equation}\label{24f14-bis}
g\in\bigoplus_{j=1}^{m}
H^{\sigma-l_j-1/2,(\sigma-l_j-1/2)/(2b);\varphi}_{+,\mathrm{loc}}(S_0,S')
\end{equation}
for
$\sigma:=p+b+n/2-2b\varkappa_k$ and some function parameter $\varphi\in\mathcal{M}$
subject to
\begin{equation}\label{24f15}
\int\limits_{1}^{\infty}\;\frac{dr}{r\varphi^2(r)}<\infty.
\end{equation}
Then the generalized derivatives
$D_{x}^{\alpha}\partial_{t}^{\beta}u_k(x,t)$ of
the component $u_k(x,t)$ of $u$ are continuous on $\Omega_0\cup\Omega'$ whenever $0\leq|\alpha|+2b\beta\leq p$.
\end{proposition}

Here, $\Omega_0$, $\Omega'$, $S_0$, and $S'$ are the sets defined at the end of Section~\ref{sec3}.

Note the cited paper \cite[Section~4]{Los17UMJ3} defines
the number $\sigma_0$ in another way; namely, $\sigma_0$ is defined to be the smallest integer such that $\sigma_0\geq\max\{0,l_1+1,\dots,l_m+1\}$ and ${\sigma_0}/{(2b)}\in\mathbb{Z}$. However, the results of this paper remain true without the assumption ${\sigma_0}/{(2b)}\in\mathbb{Z}$. To make sure of this, it is enough to use Zhitarashu and Eidelman's result \cite[Theorem~5.7]{ZhitarashuEidelman98} instead of Solonnikov's theorem in the proof of \cite[Theorem~1]{Los17UMJ3}.

Let us turn to the proof of Theorem~\ref{24th4.1}. Choose an integer $k\in\{1,\dots,N\}$ arbitrarily. To show that $u_k$ satisfies condition (a) of Definition~\ref{main-def}, we use the inclusion
$$
f\in(H^{\sigma_1,\sigma_1/(2b);\varphi_1}_{+,\mathrm{loc}}
(\Omega,\emptyset))^N
$$
from hypothesis~\eqref{24f12} of this theorem and apply Proposition~\ref{24prop1} in the case where $\Omega_0=\Omega$, $\Omega'=S_0=S'=\emptyset$, $p=2b\varkappa_k$, and $\varphi=\varphi_1$.
Then $\sigma_1:=b+n/2=\sigma$, and condition~\eqref{p-condition} is satisfied because $\sigma_3:=-b+n/2>\sigma_0$ by a hypothesis of Theorem~\ref{24th4.1}. Now we conclude by Proposition~\ref{24prop1} that $u_k$ satisfies condition~(a).

To demonstrate that $u_k$ satisfies condition (b), we make use of the hypotheses
$$
f\in(H^{\sigma_2,\sigma_2/(2b);\varphi_2}_{+,\mathrm{loc}}
(S_{\varepsilon},S))^N
$$
and \eqref{24f12bis} of Theorem~\ref{24th4.1} and hence apply Proposition~\ref{24prop1} in the case where $\Omega_0=S_\varepsilon$, $\Omega'=S_0=S$, $S'=\emptyset$, $p=l_0+2b\varkappa_k$, and $\varphi=\varphi_2$. Then $\sigma_2:=l_0+b+n/2=\sigma$, and condition~\eqref{p-condition} is satisfied because $\sigma_2>\sigma_0$ by a hypothesis of Theorem~\ref{24th4.1}. Thus, we conclude by Proposition~\ref{24prop1} that $u_k$ satisfies condition~(b).

Finally, to show that $u_k$ satisfies condition (c), we use the inclusion
$$
f\in(H_{+,\mathrm{loc}}^{\sigma_3,\sigma_3/(2b);\varphi_3}
(G_{\varepsilon},G))^N
$$
from hypothesis~\eqref{24f12} and apply Proposition~\ref{24prop1} in the case where $\Omega_0=G_{\varepsilon}$, $\Omega'=G$, $S_0=S'=\emptyset$,
$p=2b\varkappa_k-2b$, and $\varphi=\varphi_3$. Then $\sigma_3:=-b+n/2=\sigma$, and condition~\eqref{p-condition} is satisfied because $\sigma_3>\sigma_0$ by a hypothesis of Theorem~\ref{24th4.1}. We conclude by Proposition~\ref{24prop1} that the derivative $D^\alpha_x\partial^\beta_t u_k$ is continuous on $G_{\varepsilon}\cup G$ whenever $|\alpha|+2b\beta\leq 2b\varkappa_k-2b$.  This means in the $|\alpha|=0$ case that $u_k$ satisfies condition~(c).

Theorem~\ref{24th4.1} is proved.


\begin{thebibliography}{99}

\bibitem{ZhitarashuEidelman98}
S. D. Eidel'man, N. V. Zhitarashu,
\emph{Parabolic Boundary Value Problems},
Operator Theory: Advances and Applications, vol.~101, Birkh\"aser, Basel, 1998.

\bibitem{Hermander63}
L. H\"ormander, \emph{Linear Partial Differential Operators}, Grundlehren Math.
Wiss., Band~116, Springer, Berlin, 1963.

\bibitem{Ilin60}
V. A. Ilin, \emph{The solvability of  of mixed problems for hyperbolic and parabolic equations},
Russian Math. Surveys \textbf{15} (1960), no.~1, 85--142.

\bibitem{IlinKalashOleinik62}
A. M. Il'in, A. S. Kalashnikov, O. A. Oleinik, \emph{Linear equations of the second order of parabolic type}
Russian Math. Surveys \textbf{17} (1962), no.~3, 1--144.

\bibitem{Ivasyshen90}
S. D. Ivasyshen, \emph{Green Matrices of Parabolic Boundary-Value Problems},
Vyshcha Shkola, Kiev, 1990 (Russian).

\bibitem{Karamata30a}
J. Karamata, \emph{Sur certains "Tauberian theorems"\;de M.~M.~Hardy et Littlewood}, Mathematica (Cluj),
\textbf{3} (1930), 33--48.

\bibitem{LionsMagenes72ii}
J.-L. Lions, E. Magenes, \emph{Non-Homogeneous Boundary-Value Problems and
Applications, vol.~II}, Grundlehren Math. Wiss., Band~182, Springer, Berlin, 1972.

\bibitem{Los16JMathSci}
V. M. Los, \emph{Anisotropic H\"ormander spaces on the lateral surface of a cylinder},
J. Math. Sci.(N. Y.) \textbf{217} (2016), no.~4, 456 -- 467.

\bibitem{Los16UMJ6}
V. M. Los, \emph{Theorems on isomorphisms for some parabolic initial-boundary-value problems
in Hormander spaces: limiting case},
Ukrainian Math. J. \textbf{68} (2016), no.~6, 894--909.

\bibitem{Los16UMJ9}
V. M. Los, \emph{Classical solutions of the
parabolic initial-boundary value problems and  H\"ormander spaces},
Ukrainian Math. J. \textbf{68} (2016), no.~9, 1229--1239.

\bibitem{Los16UMJ11}
V. M. Los, \emph{Sufficient conditions for the solutions of general parabolic initial-boundary-value problems to be classical},
Ukrainian Math. J. \textbf{68} (2017), no.~11, 1756--1766.

\bibitem{Los17UMJ3}
V. M. Los, \emph{Systems parabolic in Petrovskii's sense in H\"ormander Spaces},
Ukrainian Math. J. \textbf{69} (2017), no.~3, 426--443.

\bibitem{Los17MFAT2}
V. Los,\emph{ Initial-boundary value problems for two-dimensional parabolic equations in H\"ormander spaces},
Methods Funct. Anal. Topology, \textbf{23} (2017), no.~2, 177--191.

\bibitem{LosMikhailetsMurach17CPAA}
V. Los V., V. A. Mikhailets, A. A. Murach, \emph{An isomorphism theorem for parabolic problems
in H\"ormander spaces and its applications},
Commun. Pur. Appl. Anal, \textbf{16} (2017), no.~1, 69--97.

\bibitem{LosMikhailetsMurach19arXiv}
V. Los V., V. A. Mikhailets, A. A. Murach, \emph{Parabolic problems in generalized Sobolev spaces},
arXiv:1907.04283.

\bibitem{LosMurach13MFAT2}
V. Los, A. A. Murach, \emph{Parabolic problems and interpolation with a function
parameter}, Methods Funct. Anal. Topology, \textbf{19} (2013), no.~2, 146--160.

\bibitem{LosMurach17OpenMath}
V. Los, A. A. Murach, \emph{Isomorphism theorems for some parabolic initial-boundary
value problems in H\"ormander spaces},
Open Mathematics \textbf{15} (2017), 57--76.

\bibitem{Lunardi1995}
A. Lunardi,  \emph{Analytic semigroups and optimal regularity in parabolic problems},
Birkhauser Verlag, Basel, 1995.

\bibitem{MikhailetsMurach12BJMA2}
V. A. Mikhailets, A. A. Murach, \textit{The refined Sobolev scale,
inter\-po\-la\-tion, and elliptic problems}, Banach J. Math. Anal. \textbf{6} (2012), no.~2, 211--281.

\bibitem{MikhailetsMurach14}
V. A. Mikhailets, A. A. Murach, \emph{Ho\"rmander spaces, interpolation, and elliptic problems}, De Gruyter, Berlin, 2014.

\bibitem{Petrovskii38}
I. G. Petrovskii, \emph{On the Cauchy problem for systems of partial differential equations
in the domain of non-anallytic functions}, Bull. Mosk. Univ., Mat. Mekh., \textbf{1} (1938), no.~7, 1--72 (Russian).

\bibitem{Solonnikov65}
V. A. Solonnikov, \emph{On boundary value problems for linear parabolic systems of differential
equations of a general form}, Trudy Mat. Inst. Steklov, \textbf{83} (1965), 3–-163  (Russian).

\bibitem{VolevichPaneah65}
L. R. Volevich and B. P. Paneah, \emph{Certain spaces of generalized functions and embedding theorems} (Russian), Uspehi Mat. Nauk \textbf{20}
(1965), no.~1, 3--74 (English translation in: Russian Math. Surveys \textbf{20} (1965), no.~1, 1--73).

\end{thebibliography}
\end{document}